\documentclass[11pt, a4paper]{article}
\usepackage{amsmath,amsfonts,amssymb,amsthm}
\newtheorem{thm}{Th\'eor\`eme}

\newtheorem{prop}{Proposition}
\newtheorem{coro}{Corollaire}%[section]
\newenvironment{remark}{\smallskip\noindent{\bf Remarque. }}{\smallskip}
\newtheorem{conj}{Conjecture}
\newtheorem{lem}{Lemme}

\begin{document}
\title{Sur une g\'en\'eralisation des coefficients binomiaux}
\author{Fr\'ederic Jouhet, Bodo  Lass et Jiang Zeng\\
\small Institut Girard Desargues,
Universit\'e Claude Bernard (Lyon 1)\\[-0.8ex]
\small 43, bd du 11 Novembre 1918,
69622 Villeurbanne cedex, France\\[-0.8ex]
\small \texttt{\{jouhet, lass, zeng\}@ euler.univ-lyon1.fr}\\[-0.8ex]}
\date{\small 2000 Mathematics Subject Classification : 05A10, 33C20}

\maketitle
\centerline{\it Combinatorica lux mea}

\begin{abstract}
We prove a recent conjecture of Lassalle about positivity and
integrality of  coefficients in some polynomial expansions.
We also give a combinatorial interpretation of those numbers.
Finally, we show that this question is closely related to the fundamental 
problem
of calculating the linearization coefficients for binomial coefficients.
\end{abstract}

%%%%%%%%%%%%%%%%%%%%%%%%%%%%%%%%%%%%%%%%%%

\section{Introduction}
Une partition $\mu= (\mu_1\geq \mu_2\geq \cdots \geq \mu_l>0)$  de $n$ est
une suite d\'ecroissante d'entiers strictement positifs de somme
$n=|\mu|$. Le nombre $l=l(\mu)$ est appel\'e
la longueur de $\mu$. Pour tout $i\geq 1$, l'entier ${m}_{i} (\mu)  =
\textrm{card} \{j: {\mu }_{j}  = i\}$ est la multiplicit\'e de
$i$ dans $\mu$. D\'efinissons
\[{z}_{\mu}  = \prod\limits_{i \ge  1}^{}
{i}^{{m}_{i}(\mu)} {m}_{i}(\mu) !  .\]
Pour  $n\geq 1$ les factorielles \emph{montantes} et
\emph{descendantes} sont d\'efinies comme suit~:
$$ \langle x\rangle_n = x (x-1) \cdots (x-n+1),\qquad
(x)_n= x (x+1) \cdots (x+n-1).
$$
Notons  que $\langle -x\rangle_n=(-1)^n(x)_n$ et que les coefficients
binomiaux valent  $\binom{x}{n}  =  {\langle x\rangle}_{n}/n!$.
Dans ses travaux sur les \emph{polyn\^omes de Jack}~\cite{La02} Lassalle
a r\'ecemment pos\'e la conjecture
suivante.
\begin{conj}
Soit $X$ une ind\'etermin\'ee, $m$ et $n$ deux
 entier strictement positifs et
${\bf r}=(r_1,\ldots,r_m)$ une suite d'entiers positifs telle que 
$|{\bf r}| = \sum_{i = 1}^m r_i>0$. On a
\begin{equation}\label{eq:las}
\sum_{|\mu| = n} \frac{X^{l(\mu) - 1}} {z_{\mu }}  \left( \sum_{i
= 1}^{l(\mu )} \prod_{k = 1}^m \frac{{(\mu_i)}_{r_k}}{r_{k}!} \right) =
\frac{1}{|{\bf r}|}  \sum_{k = 1}^{\min (n,|r|)}
c_k^{({\bf r})} \binom{X +n - 1}{n - k} ,
\end{equation}
o\`u les coefficients $c_k^{({\bf r})}$ sont 
des entiers positifs \`a d\'eterminer.
\end{conj}
Remarquons d'abord que le membre de gauche de (\ref{eq:las}) est un
polyn\^ome
en $X$ de degr\'e $n-1$, donc il peut \^etre d\'evelopp\'e dans la base
$\{\binom{X+n-1}{n-k}\}$ ($1\leq k\leq n$) d'une seule fa\c con.
Ceci implique l'existence et l'unicit\'e des coefficients rationnels
$c_k^{({\bf r})}$ au membre de droite de (\ref{eq:las}).

Comme nous allons le d\'emontrer, les nombres~$c_k^{({\bf r})}$ sont 
en fait des entiers positifs et ind\'ependants de $n$.
Pour $m=1$ et $m=2$  
les coefficients $c_k^{({\bf r})}$ ont \'et\'e d\'etermin\'es 
et la conjecture a \'et\'e v\'erifi\'ee~(voir~\cite{Ei,La00,La02,Z}). 
Dans le premier cas,  
le nombre $c_k^{(r_1)}$ est  un cofficient binomial
$$
c_k^{(r_1)}=\binom{r_1}{k},
$$
et dans le deuxi\`eme cas Lassalle~\cite{La02} a obtenu 
plusieurs formules
exprimant $c_k^{(r_1, r_2)}$, qui se r\'eduisent au cas pr\'ec\'edent lorsque
$r_2=0$.  Donc les coefficients $c_k^{({\bf r})}$ sont des extensions des 
coefficients binomiaux classiques.

L'objectif de
cet article  est de donner une solution compl\`ete de ce probl\`eme, ceci par 
\emph{trois} approches distinctes utilisant des techniques compl\`etement
diff\'erentes. Plus pr\'ecis\'ement, la section 2 donne une r\'eponse 
analytique \`a la conjecture~1, ainsi que quelques identit\'es du m\^eme 
type, ceci \`a l'aide des fonctions g\'en\'eratrices multivari\'ees. Dans la 
troisi\`eme section, nous donnons une interpr\'etation
combinatoire de l'identit\'e suivante~:
\begin{equation}\label{eq:bigeq}
\sum_{|\mu| = n} \frac{n!} {z_{\mu }} X^{l(\mu) - 1}
\sum_{i = 1}^{l(\mu )} \prod_{k = 1}^m \mu_i
\binom{\mu_i+r_k-1}{r_k-1}
=\frac{\prod_{j} r_{j}}{|{\bf
r}|}  \sum_{k = 1}^{\min (n,|r|)} c_k^{({\bf r})} k!\binom{n}{k}
(X+k)_{n-k}.
\end{equation}
Bien que la th\'eorie
des esp\`eces (voir~\cite{BLL} pour une introduction) nous ait  permis
d'imaginer cette d\'emonstration, nous pensons qu'une
pr\'esentation moins \'elitiste, plus populaire, accompagn\'ee de
calculs explicites, permettra de rendre les id\'ees
encore plus accessibles. Dans la derni\`ere section, nous d\'etaillons une 
troisi\`eme d\'emonstration de
la conjecture de Lassalle qui utilise le calcul aux diff\'erences et le cas
particulier $m = 1$, dont on trouve une d\'emonstration dans  \cite{La00}.
Dans ce paragraphe, nous voyons que le probl\`eme essentiel soulev\'e par la 
conjecture
de Lassalle est le calcul de certains coefficients de lin\'earisation.
Malgr\'e l'importance fondamentale de cette question, il semble
que, jusqu'\`a pr\'esent,  les coefficients de lin\'earisation ne
furent \'etudi\'es que pour les polyn\^omes orthogonaux.
C'est pourquoi nous ajoutons un traitement combinatoire
du probl\`eme dans ce paragraphe.

Afin de rendre la lecture la plus autonome  possible 
 nous rappelons ici 
quelques formules 
fr\'equemment  utilis\'ees dans la suite. 
D'abord  la formule binomiale peut s'\'ecrire~:
\begin{equation}\label{eq:binom}
(1-x)^{-\alpha}=\sum_{n\geq 0}\frac{(\alpha)_n}{n!}x^n.
\end{equation}
Nous aurons aussi besoin de la transformation suivante, qui
est un  cas limite de la formule de Whipple \cite[p. 142]{AAR}~:
\begin{equation}\label{whip}
{}_3F_2\!\left[\begin{matrix}-n,a,b\\c,d\end{matrix};1\right]=
\frac{(c-a)_n}{(c)_n}{}_3F_2\!\left[\begin{matrix}-n,a,d-b\\d,a+1-n-c\end{matrix};1\right],
\end{equation} 
et qui se r\'eduit \`a la formule de sommation 
de Chu-Vandermonde lorsque $b=d$~:
\begin{equation}\label{cvd}
{}_2F_1\!\left[\begin{matrix}-n,a\\c\end{matrix};1\right]=
\frac{(c-a)_n}{(c)_n},
\end{equation} 
o\`u 
\[
{}_pF_q\!\left[\begin{matrix}a_1,a_2,\dots,a_p\\
b_1,b_2,\dots,b_q\end{matrix};z\right]=
\sum _{k\geq 0} \frac {(a_1)_k\dots(a_p)_k}
{(b_1)_k\dots(b_q)_k}\frac{z^k}{k!}.
\]
est la d\'efinition  des fonctions 
hyperg\'eom\'etriques classiques. 

%%%%%%%%%%%%%%%%%%%%%%%%%%%%%%%%%%%%%%%%%%
\section{Fonctions g\'en\'eratrices}
En multipliant  le membre de gauche de (\ref{eq:las})
par $t^nx_1^{r_1}\ldots x_m^{r_m}$ et en
sommant sur $n\geq 1$ et les entiers 
$r_1,\ldots, r_m\geq 0$ tels que $|{\bf r}|\neq 0$,
par la  formule
binomiale (\ref{eq:binom}),
nous sommes amen\'es \`a \'evaluer l'expression
$$
\sum_{|\mu|\geq 1} t^{|\mu|}
{\frac{{X}^{l(\mu)-1}}{{z}_{\mu }}}
\sum_{i = 1}^{l(\mu)}\left(
\prod_{l=1}^m (1-x_l)^{-\mu_i}-1\right).
$$

\begin{lem} Soit $y$ une ind\'etermin\'ee, alors
\begin{equation}\label{eq:gauche}
\sum_{|\mu|\geq 1} t^{|\mu|}
{\frac{{X}^{l(\mu)-1}}{{z}_{\mu }}}
\sum_{i = 1}^{l(\mu)}(y^{\mu_i}-1)=\sum_{n\geq 1} t^n \sum_{k=1}^n
\binom{X+n-1}{n-k}
\frac{(y - 1 )^k}{k}.
\end{equation}
\end{lem}
\begin{proof}[Preuve]  Toute partition $\mu$ non nulle correspond de fa\c con 
biunivoque \`a une suite non nulle \`a support fini ${\mathbf 
m}=(m_1,m_2,\ldots)$ telle que
  $\mu=(1^{m_1}2^{m_2}\ldots)$. On a donc 
\begin{eqnarray}
\sum_{|\mu|\geq 1} t^{|\mu|}{\frac{{X}^{l(\mu)-1}}{{z}_{\mu }}}
\sum_{i = 1}^{l(\mu)}y^{\mu_i}
&=&\sum_{\mathbf m}X^{-1}\prod_{j\geq 
1}\left(\frac{Xt^j}{j}\right)^{m_j}\frac{1}{m_j!}\sum_{i\geq 1}m_iy^{i} 
\nonumber\\
&=& \sum_{i\geq 1}y^i
\left(\sum_{m_i\geq 0}m_i\left(\frac{Xt^i}{i}\right)^{m_i}
\frac{X^{-1}}{m_i!}\right)\cdot \prod_{j\neq i}\sum_{m_j\geq 
0}\left(\frac{Xt^j}{j}\right)^{m_j}\frac{1}{m_j!}\nonumber\\
&=&\sum_{i\geq 1}\frac{(yt)^i}{i}\exp\left(\frac{Xt^i}{i}\right)
\prod_{j\neq i}\exp\left({\frac{Xt^j}{j}}\right)\nonumber\\
&=&
(1-t)^{-X}
\log ( 1 -yt)^{-1}.\label{fg}
\end{eqnarray}
Par soustraction du terme correspondant \`a $y=1$, nous obtenons
\begin{eqnarray}
\sum_{|\mu|\geq 1} t^{|\mu|}
{\frac{{X}^{l(\mu)-1}}{{z}_{\mu }}}
\sum_{i = 1}^{l(\mu)}(y^{\mu_i}-1)&=&
(1-t)^{-X}\log \left( 1 - \frac{t}{1-t}(y - 1) \right)^{-1}\nonumber\\
&=&\sum_{k\geq 1}   (1-t)^{-X-k}\frac{t^k( y- 1 )^k }{k}\nonumber\\
&=&\sum_{n\geq 1} t^n \sum_{k=1}^n
\binom{X+n-1}{n-k}
\frac{(y - 1 )^k}{k},
\label{eq:frederic2}
\end{eqnarray}
ce qui ach\`eve la d\'emonstration.
\end{proof}
Notons, pour toute fonction multivari\'ee $f$, par $[x_1^{r_1} \cdots 
x_m^{r_m}]f(x_1 \cdots x_m)$ le coefficient de $x_1^{r_1} \cdots x_m^{r_m}$ 
dans $f$. 
Nous d\'eduisons donc de (\ref{eq:frederic2}), en posant 
$y=1/(1-x_1)(1-x_2)\cdots (1-x_m)$,
 le r\'esultat  suivant.

\begin{thm}
Soient $c_k^{({\bf r})}$ les nombres rationnels d\'efinis par (\ref{eq:las}).
Alors 
\begin{equation}\label{eq:frederic3}
\frac{c_k^{({\bf r})} }{|{\bf r}|}=[x_1^{r_1} \cdots x_m^{r_m}]\;\frac{1}{k} 
\left( \frac{1}{(1-x_1) \cdots (1-x_m)} - 1 \right)^k.
\end{equation}
En particulier,   ${kc_k^{({\bf r})} }/{|{\bf r}|}$  est un entier 
positif et  ne  d\'epend pas de~$n$.
\end{thm}
Nous en d\'eduisons donc une preuve de la conjecture~1 de Lassalle.
\begin{coro}
Les nombres $c_k^{({\bf r})}$ sont des  entiers positifs.
\end{coro}

En effet, le th\'eor\`eme~1 implique que
\begin{eqnarray}
c_k^{({\bf r})}&=&[x_1^{r_1} \cdots x_m^{r_m}]\; \frac{d}{dz} \Biggl|_{z=1} 
\frac{1}{k}
\left( \frac{1}{(1-zx_1) \cdots (1-zx_m)} - 1 \right)^k
\label{expr1}\\
&=&
[x_1^{r_1} \cdots x_m^{r_m}]\;\left( \frac{1}{(1-x_1) \cdots (1-x_m)} - 1 
\right)^{k-1}
\frac{\frac{x_1}{1-x_1}+\cdots + \frac{x_m}{1-x_m}}
{(1-x_1) \cdots (1-x_m)}.\nonumber
\end{eqnarray}
La derni\`ere expression montre clairement  que
$c_k^{({\bf r})}\in\mathbb{N}$.

Il est aussi possible de d\'eduire le corollaire au moyen des
 \emph{fonctions sym\'etriques
homog\`enes} 
sur $\{x_1,\dots,x_m\}$, qui 
sont d\'efinies~\cite{JZ, Ma} par la fonction g\'en\'eratrice~:
\begin{equation*}\label{fgh}
\sum_{n\geq 0}h_n(x_1,\ldots,x_m)z^n=\prod_{i=1}^m(1-zx_i)^{-1},
\end{equation*}
et donc ceci, \`a l'aide de (\ref{expr1}), permet d'\'ecrire :
\begin{eqnarray}
\sum_{k\geq 1}\sum_{r_1,\ldots,r_m\geq 0}c_k^{({\bf r})}t^kx_1^{r_1}\cdots
x_m^{r_m}&=&-\frac{d}{dz} \Biggl|_{z=1}\log\left(1-t\sum_{n\geq
1}h_n(x_1,\ldots,x_m)z^n\right)\nonumber\\
&=&\frac{d}{dz} \Biggl|_{z=1}\;\sum_{n\geq
1}z^n\sum_{|\lambda|=n}t^{l(\lambda)}\binom{l(\lambda)-1}{m_1(\lambda),m_2(\lambda),\ldots}h_\lambda(x_1,\ldots,x_m)\nonumber\\
&=&\sum_{\lambda}t^{l(\lambda)}|\lambda|\binom{l(\lambda)-1}{m_1(\lambda),m_2(\lambda),\ldots}h_\lambda(x_1,\ldots,x_m)\label{waring1},
\end{eqnarray}
ce qui montre aussi que $c_k^{({\bf r})}\in\mathbb{N}$. Notons que le membre 
de droite de (\ref{waring1}) s'apparente au
d\'eveloppement de la $n$i\`eme fonction sym\'etrique \emph{puissance}
$p_n(x_1,\ldots,x_m)$ dans la base
des fonctions sym\'etriques homog\`enes donn\'e par la formule de Waring
\cite{JZ, Macmahon}.\\

D'autre part, en d\'eveloppant le membre de droite de 
(\ref{eq:frederic3}) par la formule binomiale, nous obtenons
\begin{eqnarray*}
&&\frac{(-1)^k}{k}+\frac{1}{k}\sum_{i\geq 1}(-1)^{k-i}
\binom{k}{i}(1-x_1)^{-i}\cdots
(1-x_m)^{-i}\\
&=&\sum_{{\bf r}\atop |{\bf r}|> 0}{\sum_{i\geq 1}}
\frac{(-1)^{k-i}}{i}\binom{k-1}{i-1}
\prod_{l=1}^m\binom{r_l+i-1}{r_l}x_l^{r_l},
\end{eqnarray*}
ce qui donne, en extrayant le coefficient de $x_1^{r_1}\ldots x_m^{r_m}$,
 le r\'esultat suivant
\begin{coro} On a la formule explicite pour $c_k^{(\bf r)}$~:
\begin{eqnarray} 
c_k^{(\bf r)}&=&|{\bf r}|\sum_{i\geq 1}\frac{(-1)^{k-i}}{i}\binom{k-1}{i-1}
\prod_{l=1}^m\binom{r_l+i-1}{r_l}\label{eq:explicit}\\
&=&\sum_{j=1}^m\sum_{i=1}^k
(-1)^{k-i}\binom{k-1}{i-1}\binom{i+r_j-1}{r_j-1}
\prod_{l=1,l\neq j}^m\binom{r_l+i-1}{r_l}.\label{eq:entiere}
\end{eqnarray}
\end{coro}

En particulier, pour $m=1$ et $m=2$, la formule~(\ref{eq:explicit}) permet de 
retrouver les  deux expressions 
explicites de Lassalle~\cite{La02}. En fait, pour $m=1$ 
la formule 
(\ref{eq:frederic3}) se r\'eduit directement  \`a
\begin{equation}\label{eq:m=1}
\frac{k}{r_1}c_k^{(r_1)}=[x_1^{r_1}]\;x_1^k(1-x_1)^{-k}=[x_1^{r_1}]\;\sum_{l\geq 
k}\binom{l-1}{k-1}x_1^l\Longrightarrow c_k^{(r_1)}=\binom{r_1}{k}.
\end{equation}
Pour $m=2$ la formule~(\ref{eq:explicit}) s\'ecrit
 \begin{eqnarray*}
c_k^{(r_1,r_2)}&=&\frac{r_1+r_2}{k}\sum_{i=1}^k
(-1)^{k-i}\binom{k}{i}\binom{i+r_1-1}{r_1}\binom{i+r_2-1}{r_2}\\
&=&(-1)^{k-1}(r_1+r_2)\;{}_3F_2\!\left[\begin{matrix}
-k+1,r_1+1,r_2+1\\2,1\end{matrix};1\right].
\end{eqnarray*}
Appliquons deux fois  la formule~(\ref{whip})
 \`a l'expression ci-dessus, ce qui donne bien
$$c_k^{(r_1,r_2)}=\binom{r_1+r_2}{k}\;{}_3F_2\!\left[\begin{matrix}-k+1,-r_1,-r_2\\1-r_1-r_2,1\end{matrix};1\right].$$
Remarquons qu'en appliquant une troisi\`eme fois (\ref{whip}), on retrouve 
une autre expression de  \cite{La02}~:
\begin {eqnarray*}
c_k^{(r_1,r_2)}&=& \binom{r_1+r_2}{k}\binom{r_1+r_2}{r_1}
{}_3F_2\!\left[\begin{matrix}-r_1,-r_2, k-r_1-r_2\\
1-r_1-r_2, -r_1-r_2\end{matrix};1\right]\\
&=&\sum_{i \geq 0} (-1)^i
\binom{r_1+r_2-i}{k}\frac{r_1+r_2}{r_1+r_2-i}\binom{r_1+r_2-i}{i} 
\binom{r_1+r_2-2i}{r_1-i}.
\end {eqnarray*}

\begin{remark}
Lorsque tous les $r_i$ sont nuls, 
le  membre de droite  de (\ref{eq:las}) n'a pas de sens. Or
 il r\'esulte de (\ref{eq:frederic2}) avec $y=0$ que
$$
(1-t)^{-X}\log(1-t)^{-1}=\sum_{n\geq 1} t^n \sum_{k=1}^n
\binom{X+n-1}{n-k}
\frac{(- 1 )^{k-1}}{k},
$$
ce qui donne le prolongement suivant de (\ref{eq:las}) pour ${\bf r}=0$~:
\begin{equation}\label{eq;mac}
\sum_{|\mu| = n} \frac{X^{l(\mu) - 1}} {z_{\mu }}{l(\mu )}
 =\sum_{k = 1}^n\frac{(-1)^{k-1}}{k} \binom{X +n - 1}{n - k}.
\end{equation}
Cette formule est en fait la \emph{d\'eriv\'ee} d'une formule de
Macdonald~\cite[p. 26]{Ma}~:
$$
\sum_{|\mu|=n}\frac{X^{l(\mu)}}{z_\mu}=\binom{X+n-1}{n}.
$$
\end{remark}

Enfin,  en multipliant  le membre de gauche de (\ref{eq:las})
par $t^nx_1^{r_1}\ldots x_m^{r_m}$ et en
sommant sur $n\geq 1$ et les entiers 
$r_1,\ldots, r_m\geq 0$, nous obtenons
$$
\sum_{|\mu|\geq 1} t^{|\mu|}
{\frac{{X}^{l(\mu)-1}}{{z}_{\mu }}}
\sum_{i = 1}^{l(\mu)}\prod_{l=1}^m (1-x_l)^{-\mu_i},
$$
ce qui peut se d\'evelopper directement \`a l'aide de (\ref{fg}) comme suit~:
\begin{equation*}
\sum_{n\geq 0}t^n\frac{(X)_n}{n!}\sum_{k\geq 1}\frac{1}{k} \left(
\frac{t}{(1-x_1) \cdots (1-x_m)}\right)^k=\sum_{n,k\geq
1}\frac{t^n}{k}\binom{X +n-k- 1}{n - k}\prod_{l=1}^m\sum_{r_l\geq
0}\frac{(k)_{r_l}}{r_l!}x_l^{r_l},
\end{equation*}
et donc nous obtenons l'identit\'e
\begin{equation}\label{eq:las0'}
\sum_{|\mu| = n} \frac{X^{l(\mu) - 1}} {z_{\mu }}
 \sum_{i= 1}^{l(\mu )} \prod_{k = 1}^m
\frac{{(\mu_i)}_{r_k}}{r_{k}!} 
=\sum_{k = 1}^{n}\frac{1}{k}
\prod_{l=1}^m\binom{r_l+k- 1}{r_l}\binom{X +n -k- 1}{n - k}.
\end{equation}
Il est possible d'\'etablir une extension de (\ref{eq:las0'}), \`a l'aide du
coefficient $\left\langle \begin{matrix}\mu\\p\end{matrix}\right\rangle$ 
introduit par Lassalle dans \cite{La00}, et qui compte, pour
toute partition $\mu$ et tout $p\in \mathbb{N}$, le nombre de fa\c cons de
choisir $p$ \'el\'ements dans le diagramme de Ferrers de $\mu$, dont au moins
un par ligne. 
\begin{prop}
\begin{eqnarray}
&&\sum_{|\mu| = n} \left\langle\begin{matrix}\mu\\ p\end{matrix}\right\rangle
\frac{X^{l(\mu) -1}} {z_{\mu }}\left( \sum_{i= 1}^{l(\mu )} \prod_{k =
1}^m\frac{{(\mu_i)}_{r_k}}{r_{k}!} \right)\nonumber\\
&&=\sum_{k = 1}^{\min (p,|r|)}\frac{1}{k}
\left(\sum_{j=
k}^{n-p+k}\binom{j-1}{k-1}\binom{n-j-1}{p-k-1}\prod_{l=1}^m\binom{r_l+j-
1}{r_l}\right)\binom{X +p -k- 1}{p - k}.\label{eq:las0''}
\end{eqnarray}
\end{prop}
\begin{proof}[Preuve]
On a la fonction g\'en\'eratrice suivante \cite{JZ} :
$$
\sum_{p\geq 1}\left\langle\begin{matrix}\mu\\p\end{matrix}\right\rangle
 x^p=\prod_{k\geq 1}\left((1+x)^k-1\right)^{m_k(\mu)}.
 $$
Nous pouvons ainsi, comme pour (\ref{eq:las0'}),
calculer la fonction g\'en\'eratrice du membre de gauche de 
(\ref{eq:las0''}), en le multipliant par $t^nx^px_1^{r_1}\ldots x_m^{r_m}$ et 
en
sommant sur $n,p\geq 1$ et $r_1,\ldots, r_m\geq 0$ :
\begin{eqnarray*}
&&\sum_{|\mu|\geq 1} t^{|\mu|}{\frac{{X}^{l(\mu )-1}}{{z}_{\mu
}}}\sum_{p\geq 1}\left\langle\begin{matrix}\mu\\ p\end{matrix}\right\rangle
x^p\sum_{i = 1}^{l(\mu)}\prod_{l=1}^m (1-x_l)^{-\mu_i}\\
&&=\left(1-\frac{tx}{1-t}\right)^{-X}\left[\log\left(1 -\frac{t}{(1-x_1)
\cdots (1-x_m)}\right)- \log\left(1 -\frac{t(1+x)}{(1-x_1) \cdots
(1-x_m)}\right)\right].
\end{eqnarray*}
D\'eveloppons alors cette derni\`ere expression, ce qui donne :
\begin{eqnarray*}
&&\sum_{p\geq 0}\left(\frac{tx}{1-t}\right)^p\frac{(X)_p}{p!}\sum_{j\geq
1}\frac{1}{j} \left(\frac{t}{(1-x_1) \cdots
(1-x_m)}\right)^j\left((1+x)^j-1\right)\\&&=\sum_{j,k\geq 1}\sum_{p\geq
0}\frac{1}{j}\left(\frac{tx}{1-t}\right)^p\frac{(X)_p}{p!}\left(\frac{t}{(1-x_1)
\cdots (1-x_m)}\right)^j\binom{j}{k}x^k\\
&&=\sum_{j,k,p\geq
1}\frac{1}{j}\binom{j}{k}\binom{X+p-k-1}{p-k}x^p\left(\prod_{l=1}^m\sum_{r_l\geq
0}\frac{(j)_{r_l}}{r_l!}x_l^{r_l}\right)t^{p+j-k}(1-t)^{-p+k}.
\end{eqnarray*}
Mais en utilisant la formule binomiale sous la forme : 
$$(1-t)^{-p+k}=\sum_{n\geq 0}\frac{(p-k)_{n}}{n!}t^n,$$
  en rempla\c cant $n$ par $n-p-j+k$ et en extrayant le coefficient devant
$x^pt^nx_1^{r_1}\cdots x_m^{r_m}$, nous obtenons la fonction g\'en\'eratrice 
du membre de droite.
\end{proof}

\begin{remark}
Pour $p=n$, l'identit\'e (\ref{eq:las0''})
 donne  bien (\ref{eq:las0'}).
\end{remark}

%%%%%%%%%%%%%%%%%%%%%%%%%%%%%%%%%%%%%%%%%%%%%%%%%%%%%%%%%%%%%%%%%%%%%%%%
 \section{Interpr\'etation en  th\'eorie d'esp\`eces}
Imaginons que $m$ esp\`eces d'animaux, \`a savoir $r_1$~\^anes, 
$r_2$~belettes,
$r_3$~chevaux, $r_4$~daims, $r_5$~\'ecureuils, \dots, $r_m$~mulets
organisent une \emph{Table Ronde} sur le th\`eme \emph{Protection des
esp\`eces}.   \'Evidemment, l'\'ev\'enement a lieu \`a Montr\'eal
dans la salle de conf\'erence de l'UQAM, o\`u  il y a $n$~chaises \`a la
disposition des participants. Le comit\'e d'organisation les place
autour de plusieurs tables rondes, et rattache chaque chaise \`a
ses deux voisines \`a l'aide d'un ruban \'elastique pour bien fixer
l'ordre.
\'Etant donn\'ee une partition $\mu$ de $n$,
il est \'evident que le comit\'e
d'organisation a
$n!/z_\mu$ possibilit\'es pour placer $n$~chaises donn\'ees
autour de $l = l(\mu)$  tables rondes,
\`a savoir
$m_1(\mu)$~tables avec une seule place, $m_2(\mu)$~tables avec
deux places, $m_3(\mu)$~tables avec trois places, \dots,
$m_n(\mu)$~tables avec $n$~places.  Apr\`es avoir
termin\'e ce travail, le comit\'e contacte l'\^ane le plus \^ag\'e,
pour que celui-ci puisse en choisir une. Les $l(\mu)-1$~autres tables rondes
restent
donc libres et sont chacune r\'eserv\'ees pour l'un des $X$~autres
s\'eminaires programm\'es \`a
Montr\'eal. En fait, c'est le comit\'e d'organisation qui se
d\'ecide en faveur d'une des $X^{l(\mu)-1}$~r\'eservations
possibles.

Puisque l'\^ane n'a r\'eserv\'e qu'une seule table avec
$\mu_i$~places, il est d\'ej\`a assez clair que tout le monde ne
pourra pas venir, d'autant plus que les organisateurs n'acceptent
pas que deux repr\'esentants d'une m\^eme esp\`ece s'assoient sur
une m\^eme chaise. En effet, dans ce cas ils risqueraient de
chuchoter l'un avec l'autre tout le temps.  Il est cependant
tout \`a fait admissible et m\^eme, dans l'esprit de l'entente entre les
esp\`eces, d\'esirable, que des repr\'esentants de
diff\'erentes esp\`eces s'installent sur une m\^eme chaise. Par
ailleurs, toute tentative d'\'eviter cette cohabitation serait
perdue d'avance, parce que chaque esp\`ece tient beaucoup  \`a son
ind\'ependance~: notamment dans le choix des chaises. Les \^anes
d\'ecident donc d'envoyer l'une des $\binom{r_1}{a}$~d\'el\'egations
  possibles, o\`u,  \'evidemment, $1 \le a \le r_1$,
puisque chaque esp\`ece doit \^etre repr\'esent\'ee \`a
Montr\'eal. Apr\`es \^etre arriv\'es dans la salle de conf\'erence,
les commissaires choisissent $a$~chaises pour s'y installer dans
un ordre lin\'eaire. Autrement dit, l'\^ane le plus \^ag\'e de la
commission commence par choisir sa place, et les autres
commissaires s'installent, suivant l'\^age, l'un apr\`es l'autre sur
les places choisies \`a sa gauche jusqu'\`a ce que le plus jeune
\^ane de la commission s'assoit sur la place choisie \`a la droite
du doyen.   Pour les \^anes, il y a donc effectivement, en vertu de 
(\ref{cvd}),
$$
\sum_{a=1}^{r_1}      \binom{r_1}{a}  \binom{\mu_i}{a}   a
%=\mu_i    \sum_{a=1}^{r_1}    \binom{r_1}{r_1-a}  \binom{\mu_i-1}{a-1}
=
 \mu_i    \binom{\mu_i+r_1-1}{r_1-1}
% =r_1   \binom{\mu_i+r_1-1}{r_1} \nonumber
$$
mani\`eres diff\'erentes de prendre leurs fonctions. En r\'esum\'e nous avons 
\'etabli le r\'esultat suivant
%%%%%%%%%%%%%%%%%%%%%%%%%%%%%%%%%%%%%%%%%
\begin{lem}\label{com} Soit $F_k(r)$ le nombre de fa\c cons 
de choisir des  commissaires
 d'une esp\`ece
ayant  $r$ repr\'esentants et de les installer
autour d'une table ayant  $k$ chaises, alors
$$
F_k(r)=k\binom{k+r-1}{r-1}=r\binom{k+r-1}{k-1}.
$$
\end{lem}
\noindent Par cons\'equent, nous avons le r\'esultat explicite suivant~:
\begin{prop} Si $F_k ({\bf r})$ (resp.~$S_k ({\bf r})$) est 
le nombre de mani\`eres de choisir des commissaires de chaque esp\`ece et 
de les installer  autour d'une table avec $k$~places
(resp.~dont aucune ne doit rester vide),
alors on a
\begin{equation}\label{eq:comvide}
 F_k ({\bf r})=\prod_{i=1}^mF_{k} ({r_i})=
\prod_{l=1}^m r_l \binom{k+r_l-1}{r_l-1},
\end{equation}
et
\begin{equation}\label{eq:comnonvide}
S_k ({\bf r})=
\sum_{i=1}^k (-1)^{k-i}  \binom{k}{i}
\prod_{l=1}^m r_l  \binom{i+r_l-1}{r_l}.
\end{equation}
\end{prop}
En effet, d'apr\`es le lemme~\ref{com}
la formule (\ref{eq:comvide}) est \'evidente, et d'autre part
on a l'\'equivalence suivante~:
$$
 F_k ({\bf r})
 =
\sum_{i=1}^k  \binom{k}{i}   S_i ({\bf r})
\Longleftrightarrow 
S_k ({\bf r})=
\sum_{i=1}^k (-1)^{k-i} \binom{k}{i} F_i ({\bf r}),
$$
qui permet de  d\'eduire (\ref{eq:comnonvide}) par substitution 
de (\ref{eq:comvide}).

Voil\`a pourquoi le nombre total de sc\'enarios
diff\'erents  est \'egal \`a
$$
 \sum_{|\mu|=n}\frac{n!}{z_\mu}
X^{l(\mu)-1} \sum_{i=1}^{l(\mu)}  F_{\mu_i}
({\bf r}).
 $$

La r\'eunion commence bien \`a l'heure. H\'elas, c'est un
\'ecureuil qui cause les premiers retards en critiquant la
politique des chaises vides. Apr\`es un vote, les d\'el\'egations
d\'ecident donc que chaque commissaire qui trouve une chaise vide
\`a sa droite \'eloigne celle-ci et coupe le ruban \'elastique. Les
organisateurs sont dans tous leurs \'etats apr\`es avoir appris que
des rubans furent coup\'es. Mais au bout du compte, ils prennent
leur parti de la situation et arrangent les chaises vides en
plusieurs \emph{queues} derri\`ere les chaises occup\'ees tout en
laissant les autres rubans intacts.

Une belette pense qu'il aurait \'et\'e plus simple de choisir les
$k$~chaises r\'eellement occup\'ees d'abord et de les placer
autour d'une table ronde plus petite. Le doyen des daims ajoute
que cela aurait \'et\'e possible de $(k-1)!\binom{n}{k}$~mani\`eres
diff\'erentes, mais ensuite il aurait \'et\'e fort
difficile de partitionner les autres $n-k$~chaises en $c$~ordres
cycliques et $l$~ordres lin\'eaires, d'autant plus que cette
configuration devrait \^etre compt\'ee avec un facteur
$X^c \langle k\rangle_l$, puisque il faut bien placer
chaque \emph{queue} de chaises derri\`ere une des
$k$~places occup\'ees. Rien de plus simple que cela,
s'exclame un petit \'ecureuil : on obtient le polyn\^ome
de recouvrement~$C!(K_{n-k},X,k)$ du graphe orient\'e
complet~$K_{n-k}$ (\cite{CG})~! Et, gr\^ace au  th\'eor\`eme de
dualit\'e (\cite{C},~\cite{L}), le r\'esultat est \'egal \`a
\begin{equation*}
C!(K_{n-k},X,k) = (-1)^{n-k} C!(\overline{K_{n-k}},X,-X-k)  =
(-1)^{n-k} \langle-X-k\rangle_{n-k} =  (X+k)_{n-k},
\end{equation*}
puisque le graphe orient\'e~$\overline{K_{n-k}}$ sans aucun
arc n'admet qu'une seule partition en \hbox{$n-k$}~ordres
lin\'eaires. Un cheval trouve que ce n'est vraiment pas la peine
de renvoyer les lecteurs aux \oe uvres de ces jeunes vauriens
quand, en r\'ealit\'e, on utilise des r\'esultats classiques imagin\'es
par des ma\^\i tres tels Berge~(\cite{B}), Foata et
Strehl~(\cite{FS})~:
\begin{equation*}
\sum_{f: [n-k] \to [n]}  X^{{\mathrm{cyc}}\, f}
=
(X+k)_{n-k},
\end{equation*}
o\`u  la somme porte sur toutes les injections
$f: \{1,\dots,n-k\} \to \{1,\dots,n\}$ ($\mathrm{cyc}\, f$
est le nombre de cycles de~$f$).
Un grand mulet, cependant, pense qu'il serait souhaitable
de pr\'esenter une d\'emonstration \`a la lumi\`ere de
la th\'eorie des esp\`eces~:
\begin{equation*}
\exp \left[ X \sum_{i\geq 1}
(i-1)! \frac{t^i}{i!} \right]
\left[ 1+ \sum_{i\geq 1}
i! \frac{t^i}{i!}  \right]^k
=
\exp\left[ -X\log(1-t) \right]  \left[ 1-t \right]^{-k}
 \nonumber
\end{equation*}
\begin{equation*}
 =   \left[ 1-t \right]^{-X-k}
 =
 1 + \sum_{i\geq 1}
(X+k)_i \frac{t^i}{i!}.
\nonumber
\end{equation*}

Apr\`es une halte contemplative, un \^ane remarque que l'on
aurait, par ailleurs, \'etabli deux formules nouvelles pour le
nombre de sc\'enarios diff\'erents~:  l'une, plus difficile, correspondant 
\`a ce qui vient d'\^etre discut\'e, et l'autre, plus simple,
correspondant au cas o\`u l'on n'aurait pas coup\'e de ruban \'elastique.

\begin{thm} Le nombre total de sc\'enarios diff\'erents peut 
s'exprimer  comme suit~:
\begin{eqnarray}
\sum_{|\mu|=n}\frac{n!}{z_\mu}
X^{l(\mu)-1} \left( \sum_{i=1}^{l(\mu)}
F_{\mu_i} ({\bf r})  \right)
&=&
\sum_{k=1}^n     F_k ({\bf r})
(k-1)! \binom{n}{k} (X)_{n-k}
\label{eq:vraif} \\
&=&
\sum_{k=1}^n   S_k ({\bf r}) (k-1)!\binom{n}{k}
(X+k)_{n-k}.  \label{eq:vrais}
\end{eqnarray}
\end{thm}
Les identit\'es (\ref{eq:vrais}) et (\ref{eq:vraif}) correspondent
respectivement aux identit\'es (\ref{eq:bigeq}) et (\ref{eq:las0'}).\\
On en d\'eduit alors que 
\begin{equation} \label{eq:last}
c_k^{({\bf r})}
 =
 \frac{|{\bf r}|} {k\cdot \prod_{j} r_{j}}  S_k ({\bf r})
 =
 \sum_{j=1}^m \frac{S_k ({\bf r})\cdot r_j}
 {k \cdot r_1 \cdots r_m},
\end{equation}
ce qui montre que $c_k^{({\bf r})}$ est positif et  
ne  d\'epend pas de~$n$, 
et par substitution de (\ref{eq:comnonvide}), on retrouve les formules du 
corollaire~2,  dont la derni\`ere, \`a savoir (\ref{eq:entiere}),
 montre que $c_k^{({\bf r})}$ est un entier.

En fait, nous pouvons renforcer le dernier r\'esultat, c'est-\`a-dire la 
conjecture de Lassalle.
Supposons que les $m$ esp\`eces d'animaux 
 soient num\'erot\'ees de 1 \`a $m$. 
La 1\`ere esp\`ece est donc celle des \^anes et 
nous pouvons parler de la $j$i\`eme esp\`ece avec $1\leq j\leq m$.
\begin{thm} Etant donn\'ees $m$ esp\`eces ayant respectivement 
$r_1,\ldots, r_m$ repr\'esentants et une table entour\'ee de
$k$ chaises num\'erot\'ees de 1 \`a $k$,
le nombre de fa\c  cons de choisir des commissaires de chaque esp\`ece et
de les installer autour de la  table ayant $k$ chaises
de sorte qu'aucune chaise ne soit vide et que le doyen de la
commission de la $j$i\`eme esp\`ece soit 
install\'e sur la chaise num\'ero $k$, et que le doyen  
de toute autre esp\`ece fasse partie de sa propre commision est donn\'e par
$$
T_k({\bf r};j)=\frac{S_k ({\bf r})\cdot r_j}{k \cdot r_1 \cdots r_m}, \qquad 
\hbox{pour}\quad 1\leq j\leq m.
$$
\end{thm}
\begin{proof}[Preuve] Evidemment nous pouvons supposer sans perdre de 
g\'en\'eralit\'e
que $j=1$. Il s'agit donc de d\'emontrer que
$$
S_k ({\bf r})=kr_2\cdots r_mT_k({\bf r};1).
$$
Le doyen de la commission des \^anes (ce n'est pas forc\'ement le
doyen de tous les $r_1$~\^anes!) 
 choisit, parmi toutes les $k$~chaises (et pas seulement parmi
les chaises occup\'ees par les \^anes), celle qui porte le plus grand 
num\'ero (ici le num\'ero $k$)  pour pr\'esider la s\'eance. Les autres 
esp\`eces,
cependant, sont oblig\'ees d'installer le doyen de toute leur
esp\`ece (et pas seulement le doyen de leur commission) sur la
chaise la plus grande parmi celles occup\'ees par des commissaires
de leur esp\`ece (ce n'est pas forc\'ement la plus grande de
toutes les  $k$~chaises!).\end{proof}

La formule  (\ref{eq:last}) peut s'interpr\'eter comme suit.
Les commissaires des autres esp\`eces approuvent le principe de
pr\'esidence sugg\'er\'e par le  doyen de la commission des
\^anes, mais, naturellement, ils insistent sur l'id\'eal de
l'\'egalit\'e de toutes les esp\`eces. Suivant l'exemple de l'UE,
on se d\'ecide donc en faveur d'une pratique du tourniquet, ce qui
augmente le nombre total de sc\'enarios diff\'erents \`a
$c_k^{({\bf r})}=\sum_{j=1}^mT_k({\bf r};j)$
pour chaque $k \in {\mathbb N}$. 

Tout le
monde est enchant\'e; seul le jeune mulet revient sur sa question, \`a savoir 
comment on pourrait installer les commissaires de fa\c con
surjective tout en respectant l'ind\'ependance de toutes les
esp\`eces. Un \'ecureuil pense que l'on pourrait utiliser la
th\'eorie des esp\`eces virtuelles (\cite{BLL}, sect.~2.5.) pour
r\'esoudre ce probl\`eme difficile~\dots

Avant la deuxi\`eme conf\'erence, des militants antimondialisation
se sont infiltr\'es dans la salle de conf\'erence pour mettre des
graffitis sur $p$~chaises, et notamment sur au moins une chaise \`a
chaque table ronde. Ceci augmente le nombre de sc\'enarios \`a
$$
\sum_{\left|{\mu }\right| = n}
\left\langle
\begin{matrix}\mu\\p\end{matrix}
\right\rangle
{\frac{n!}{{z}_{\mu }}} {X}^{l(\mu ) - 1}
\left(\sum_{i = 1}^{l(\mu)}   F_{\mu_i} ({\bf r})
\right).
$$
Le jeune mulet sugg\`ere qu'il faudrait commencer par
choisir les $p$~chaises d\'egrad\'ees, et noter~$i$ (resp.~$j$) le
nombre de chaises d\'egrad\'ees (resp.~non endommag\'ees) parmi
les $\mu_i$~chaises autour de la table choisie par le doyen des
\^anes pour la conf\'erence.   Si l'on \'eloigne les chaises non
endommag\'ees de toutes les autres tables, alors il y a
$$
\binom{n}{p} \sum_{i=1}^p \sum_{j=0}^{n-p} F_{i+j}
({\bf r}) (i+j-1)!\binom{p}{i}\binom{n-p}{ j} (X)_{p-i}
$$
possibilit\'es diff\'erentes. Comme il y a
$(p-i)(p-i+1)(p-i+2)\cdots (n-i-j-1) = (p-i)_{n-p-j}$ mani\`eres
diff\'erentes de r\'eintroduire les $n-p-j$~chaises, on a d\'emontr\'e 
l'identit\'e suivante
\begin{equation*}
 \sum_{\left|{\mu }\right| = n}
 \left\langle
\begin{matrix}\mu\\p\end{matrix}
\right\rangle
{\frac{n!}{{z}_{\mu }}} {X}^{l(\mu ) - 1}
\left(\sum_{i = 1}^{l(\mu)}   F_{\mu_i} ({\bf r})
\right)
\nonumber
\end{equation*}
\begin{equation*}\label{eq:graf}
\, = \, \binom{n}{p} \sum_{i=1}^p \sum_{j=0}^{n-p}
F_{i+j} ({\bf r}) (p-i)_{n-p-j} (i+j-1)!
\binom{p}{i}\binom{n-p}{j} (X)_{p-i},
\end{equation*}
qui est exactement l'identit\'e (\ref{eq:las0''}).

%%%%%%%%%%%%%%%%%%%%%%%%%%%%%%%%%%%%%%%%%%%%%%%%%%%%%%%%%%%%%%%
\section{Liens avec les coefficients de lin\'earisations}
Remarquons d'abord qu'en posant  $X=0$
dans l'\'equation (\ref{eq:las}) nous obtenons 
\begin{equation}\label{eq:poly}
\prod_{i=1}^m\frac{(n)_{r_i}}{r_j!} = \frac{1}{|{\bf
r}|}\sum_{k=0}^{|{\bf r}|} k\,c_k^{({\bf r})}\,\frac{\langle
n\rangle_{k}}{k!}.
\end{equation}
Comme $c_k^{({\bf r})}$ est ind\'ependant de $n$, 
la d\'etermination de 
$c_k^{({\bf r})}$ appara\^{\i}t donc comme le calcul des 
coefficients de d\'eveloppement du polyn\^ome $(x)_{r_1}\ldots (x)_{r_m}$
dans la base  $(\langle x\rangle_{k})_{k\geq 0}$. De plus, 
si nous pouvons d\'emontrer autrement que les nombres 
$c_k^{({\bf r})}$ sont  ind\'ependants de $n$, cette approche  fournirait une 
nouvelle
preuve de la conjecture de Lassalle.

Comme dans le paragraphe pr\'ec\'edent, nous
consid\'erons $r_1$~\^anes, $r_2$~belettes, \dots, $r_m$~mulets,
qui veulent s'asseoir sur $x$~chaises. De nouveau, deux
repr\'esentants d'une m\^eme esp\`ece ne sont pas autoris\'es
\`a choisir la m\^eme chaise. Il est cependant admissible que
des repr\'esentants de  diff\'erentes esp\`eces s'installent
sur une m\^eme chaise. En fait, ceci est, en g\'en\'eral,
m\^eme in\'evitable puisque les esp\`eces sont
ind\'ependantes dans leur choix des chaises.
Voil\`a pourquoi le nombre de sc\'enarios
possibles est \'egal \`a $\langle x \rangle_{r_1}
\langle x \rangle_{r_2} \cdots \langle x \rangle_{r_m}$.

Soit $E = [r_1] \uplus [r_2] \uplus \dots \uplus [r_m]$
l'union disjointe des repr\'esentants de toutes les esp\`eces.
Appelons un sous-ensemble $T \subseteq E$ \emph{transversal} 
si $\textrm{card} (T \cap [r_i]) \in \{0,1\}$
pour tout $i \in \{1,2,\dots,m\}$. Les transversaux
de~$E$ sont \'evidemment les sous-ensembles
qu'on peut installer sur une seule chaise.  Ceci
d\'emontre le th\'eor\`eme suivant.
\begin{thm} 
Soit $d_k(r_1,\dots,r_m)$  le nombre de mani\`eres
diff\'erentes de partitionner $E$ en $k$~transversaux
non-vides, alors
\begin{equation}\label{eq:linm}
\langle x\rangle_{r_1}\cdots \langle x\rangle_{r_m}
=
\sum_{k\geq 0}d_k({\bf r}) \langle x\rangle_{k},
\end{equation}
En particulier,  nous avons la formule de lin\'earisation classique~:
\begin{equation}
\langle x\rangle_{r_1}\langle x\rangle_{r_2}
=\sum_{k\geq 0}\binom{{r_1}}{k}\binom{{r_2}}{k}k! \langle 
x\rangle_{{r_1}+{r_2}-k}
\label{eq:lin2}
\end{equation}
\end{thm}
\noindent
En effet, pour $m=2$, s'il y a
$k$~transversaux de cardinal deux et si le nombre total de transversaux vaut 
$r_1+r_2-k$, alors nous pouvons les choisir
de $\binom{r_1}{k} \binom{r_2}{k} k!$  fa\c cons distinctes,
 c'est-\`a-dire
$$
d_{r_1+r_2-k} (r_1,r_2) =
\binom{r_1}{k} \binom{r_2}{k} k!.
$$

Il est encore plus simple de choisir directement, de fa\c con
ind\'ependante, $m$~sous-ensembles de~$[x]$ de cardinaux
$r_1$, \dots, $r_m$, respectivement.  Ceci est possible de
$\binom{x}{r_1} \cdots \binom{x}{r_m}$ mani\`eres
distinctes et montre le th\'eor\`eme suivant.
\begin{thm} Soit
$\tilde{d}_k({\bf r})$   le nombre de mani\`eres
diff\'erentes de choisir  $m$~sous-ensembles de~$[k]$ de cardinaux
$r_1$, \dots, $r_m$, respectivement, de sorte que chaque \'el\'ement
de~$[k]$ soit choisi au moins une fois. Alors
\begin{equation}\label{eq:linbin}
\binom{x}{r_1}\cdots \binom{x}{r_m}
=
\sum_{k\geq 0} \tilde{d}_k({\bf r}) \binom{x}{k},  \qquad
\tilde{d}_k({\bf r})
=
\frac{k! \, d_k({\bf r}) } {r_1! \cdots r_m!},
\end{equation}
 En particulier, on a
$$
\tilde{d}_{r_1+r_2-k} (r_1,r_2) =
\binom{r_1+r_2-k}{k,r_1-k,r_2-k}.
$$
\end{thm}
\noindent
On peut aussi donner une preuve directe de ce dernier r\'esultat.
En effet, choisir deux  sous-ensembles $E_1$ et $E_2$ 
de~$[x]$ tels que $|E_1|=r_1$, $|E_2|=r_2$ et  $|E_1\cap E_2|=k$ 
\'equivaut \`a choisir un sous-ensemble de~$[x]$ de
cardinal~$r_1+r_2-k$ et puis le partitionner en trois blocs de
cardinaux $k$, $r_1-k$, $r_2-k$, respectivement. D'o\`u
$\tilde{d}_{r_1+r_2-k} (r_1,r_2) =
\binom{r_1+r_2-k}{k,r_1-k,r_2-k}$.

Au lieu de choisir, de fa\c con ind\'ependante,  $r_1$, \dots, $r_m$
\'el\'ements de~$[x]$ sans r\'ep\'etition, choisissons-les maintenant
avec des r\'ep\'etitions possibles.  
Comme le nombre de fa\c cons de choisir $n$ \'el\'ements dans 
$[x]$ avec des r\'ep\'etitions possibles est 
$$
\left(\binom{x}{n}\right) = \binom{x+n-1}{n}=\frac{(x)_{n}}{n!},
$$
le nombre de sc\'enarios
distincts est donc \'egal \`a $\bigl(\binom{x}{r_1}\bigr) \cdots 
\bigl(\binom{x}{r_m}\bigr)$.
Une comparaison avec (\ref{eq:poly}) montre le th\'eor\`eme suivant.
\begin{thm} Soit $\tilde{c}_k({\bf r})$  le nombre de mani\`eres
diff\'erentes de choisir   $r_1$, \dots, $r_m$   \'el\'ements de~$[k]$
avec des r\'ep\'etitions possibles, de sorte que chaque \'el\'ement
de~$[k]$ soit choisi au moins une fois, alors
\begin{equation}\label{eq:linlas}
\biggl(\binom{x}{r_1}\biggr) \cdots \biggl( \binom{x}{r_m}\biggr)
=
\sum_{k\geq 0} \tilde{c}_k({\bf r}) \binom{x}{k}. 
\end{equation}
  En particulier on a
\begin{equation}
\tilde{c}_k({r_1, r_2})=\sum_{l+k_1+k_2=k}\binom{k}{l,k_1-l,k_2-l}
\binom{r_1-1}{k_1-1}\binom{r_2-1}{k_2-1}.
\end{equation}
\end{thm}
\noindent
Il est \'evident que (\ref{eq:linlas}) et (\ref{eq:frederic3})
fournissent exactement les m\^emes interpr\'etations combinatoires
pour les nombres $c_k^{({\bf r})}$ introduits par Lassalle.

Notons que l'identit\'e (\ref{eq:lin2}) s'\'ecrit  encore
\begin{equation}\label{eq:binom2}
\frac{(x)_{r_1}}{r_1!} \frac{(x)_{r_2}}{r_2!}
=
\sum_{l\geq 0}(-1)^l  \binom{r_1+r_2-l}{l,r_1-l,r_2-l}
\frac{(x)_{r_1+r_2-l}}{(r_1+r_2-l)!}.
\end{equation}
En utilisant (\ref{eq:binom2}) dans (\ref{eq:las}) nous d\'eduisons
le r\'esultat suivant~:
\begin{lem} Les coefficients $c_k^{({\bf r})}$ satisfont
 la relation de r\'ecurrence suivante~:
\begin{equation}\label{eq:recurr}
\frac{c_k^{( r_1,r_2,r_3,\dots,r_m)}} {r_1+r_2+r_3+\cdots+r_m}
=
\sum_{l \ge 0} (-1)^l  \binom{r_1+r_2-l}{l,r_1-l,r_2-l}
\frac{c_k^{( r_1+r_2-l,r_3,\dots,r_m)}} {r_1+r_2 -l+r_3+\cdots+r_m}.
\end{equation}
En particulier, 
comme $c_k^{(r_1)}=\binom{r_1}{k}$ (voir (\ref{eq:m=1})),
les coefficients  $c_k^{({\bf r})}$   sont ind\'ependants de
$n$.
\end{lem}

Par comparaison de (\ref{eq:poly}) et (\ref{eq:linlas}) il en r\'esulte  que
$$
\tilde{c}_k({\bf r})
={ k\,c_k^{({\bf r})} }  /  {|{\bf r}|}.
$$

En vue de d\'eduire une nouvelle preuve de la conjecture de Lassalle, nous 
introduisons quelques notations suppl\'ementaires.
Pour tout polyn\^ome $P(x)$ d\'efinissons les op\'erateurs 
$E$, $I$ et $\Delta$ comme suit~:
$$E P(x) = P(x+1), \quad I P(x)=P(x)\quad \hbox{et}\quad \Delta =E-I.
$$ 
Pour tout $k\geq 0$ posons $\Delta^0(P(x))=P(x)$ et
$\Delta^{k+1}=\Delta(\Delta^k)$. 
La formule binomiale implique que
\begin{equation}\label{eq:diff}
\Delta^nP(x)
=(E-I)^nP(x)=
\sum_{k=0}^n(-1)^k\binom{n}{k}P(x+n-k),
\end{equation}
et d'autre part nous avons le d\'eveloppement de Taylor suivant~:
\begin{equation}\label{eq:taylor}
P(x)
=
\sum_{k\geq 0}\frac{\Delta^kP(0)}{ k!}\langle x\rangle_k.
\end{equation}
En vertu de la formule de Chu-Vandermonde (\ref{cvd}) on a 
\begin{equation*}
(x)_n
=
\sum_{j\geq 0}\binom{n}{j}(j)_{n-j}\langle x\rangle_j.
\end{equation*}
Ainsi
\begin{equation}\label{eq:change}
\prod_{i=1}^m(x)_{r_i}
=
\sum_{j_1,\ldots, j_m\geq 0}\prod_{i=1}^m
\binom{r_i}{j_i}(j_i)_{r_i-j_i}\langle x\rangle_{j_i}.
\end{equation}
Substituons  (\ref{eq:linm}) dans (\ref{eq:change})~:
\begin{equation}\label{eq:change2}
\prod_{i=1}^m(x)_{r_i}
=\sum_{k=0}^{|{\bf r}|}\sum_{j_1,\ldots, j_m\geq 0}\left(\prod_{i=1}^m
\binom{r_i}{j_i}(j_i)_{r_i-j_i}\right)d_k({\bf j})\langle x\rangle_{k}.
\end{equation}
D'autre part, en appliquant  directement (\ref{eq:diff}) et (\ref{eq:taylor})
avec $P(x)=(x)_{r_1}\cdots (x)_{r_m}$  nous obtenons
\begin{equation}\label{eq:connex}
\prod_{i=1}^m(x)_{r_i}
=
\sum_{k=0}^{|{\bf r}|}
\left(\frac{1}{k!}\sum_{j=0}^k(-1)^j\binom{k}{ j}
\prod_{i=1}^m (k-j)_{r_i}\right)\, \langle x\rangle_{k}.
\end{equation}
Gr\^ace au lemme~3 la comparaison de (\ref{eq:poly}) 
avec (\ref{eq:change2}) et (\ref{eq:connex})
montre le th\'eor\`eme suivant.
\begin{thm} On a d'une part
\begin{equation}\label{eq:positive}
c_k^{({\bf r})}=
\frac{r_1+\cdots +r_m}{r_1!\cdots r_m!}
\sum_{j_1,\ldots, j_m\geq 0}\left(\prod_{i=1}^m\binom{r_i}{j_i}(j_i)_{r_i-j_i}\right)
\frac{d_k({\bf j})}{k},
\end{equation}
et d'autre part la formule explicite
(\ref{eq:entiere}), c'est-\`a-dire,
\begin{equation}\label{eq:entiere'}
c_k^{({\bf r})}=\sum_{j=1}^m\sum_{i=1}^k
(-1)^{k-i}\binom{k-1}{i-1}\binom{i+r_j-1}{r_j-1}
\prod_{l=1,l\neq j}^m\binom{r_l+i-1}{r_l}.
\end{equation}
\end{thm}
\noindent Il r\'esulte   respectivement de 
(\ref{eq:positive}) et (\ref{eq:entiere'}) que
$c_k^{({\bf r})}$ est entier et positif.

\begin{remark}
Un $q$-analogue des r\'esultats de cette derni\`ere section sera trait\'e
dans un article ult\'erieur.

\end{remark}

\end{document}